\documentclass[draftclsnofoot,11pt,onecolumn]{IEEEtran}

\usepackage{latexsym}
\usepackage{amsmath}
\usepackage{amsfonts}
\usepackage{amssymb}
\usepackage{amsthm}
\usepackage{cite}
\usepackage{graphicx} 
\usepackage{subfigure} 
\usepackage{color}
\usepackage{algorithm} 
\usepackage{algorithmic} 
\usepackage{bm}
\usepackage{longtable}
\usepackage{hyperref}
\usepackage{graphicx}  

\newtheorem{theorem}{Theorem}
\newtheorem{lemma}{Lemma} 
\newtheorem{remark}{Remark} 
\newtheorem{corollary}{Corollary} 
\ifCLASSINFOpdf
\else
\fi
\begin{document}
%
\title{Coherence-Based Performance Guarantee of Regularized $\ell_{1}$-Norm Minimization and Beyond}
%
%
%

\author{

\thanks{
This work was partially supported by the National Science Foundation of China (Grant Nos. 61273020, 61673015) and the China Postdoctoral Science Foundation (Grant No. 2018M643390). (\emph{Corresponding author: Jianjun Wang.})
}

Wendong~Wang,
Feng~Zhang,
Zhi~Wang,
and~Jianjun~Wang

\thanks{W. D. Wang, F. Zhang, and J. J. Wang are with the School of Mathematics and Statistics, Southwest University, Chongqing 400715, China (e-mail: wdwang@swu.edu.cn; zhangf@email.swu.edu.cn; wjj@swu.edu.cn).}
\thanks{Z. Wang is with the School of Mathematics and Statistics and the College of Computer \& Information Science, Southwest University, Chongqing, 400715, China (e-mail: chiw@swu.edu.cn).}

}

%
%

\markboth{Preprint submitted to Journal of LATEX Templates}
{WANG \MakeLowercase{\textit{et al.}}: Coherence-Based Performance Guarantee of Regularized $\ell_{1}$-Norm Minimization and Beyond}
%



\maketitle

\begin{abstract}
In this paper, we consider recovering the signal $\bm{x}\in\mathbb{R}^{n}$ from its few noisy measurements $\bm{b}=A\bm{x}+\bm{z}$, where $A\in\mathbb{R}^{m\times n}$ with $m\ll n$ is the measurement matrix, and $\bm{z}\in\mathbb{R}^{m}$ is the measurement noise/error. We first establish a coherence-based performance guarantee for a regularized $\ell_{1}$-norm minimization model to recover such signals $\bm{x}$ in the presence of the $\ell_{2}$-norm bounded noise, i.e., $\|\bm{z}\|_{2}\leq\epsilon$, and then extend these theoretical results to guarantee the robust recovery of the signals corrupted with the Dantzig Selector (DS) type noise, i.e., $\|A^{T}\bm{z}\|_{\infty}\leq\epsilon$, and the structured block-sparse signal recovery in the presence of the bounded noise. To the best of our knowledge, we first extend nontrivially the sharp uniform recovery condition derived by Cai, Wang and Xu (2010) for the constrained $\ell_{1}$-norm minimization model, which takes the form of
\begin{align*}
\mu<\frac{1}{2k-1},
\end{align*}
where $\mu$ is defined as the (mutual) coherence of $A$, to two unconstrained regularized $\ell_{1}$-norm minimization models to guarantee the robust recovery of any signals (not necessary to be $k$-sparse) under the $\ell_{2}$-norm bounded noise and the DS type noise settings, respectively. Besides, a uniform recovery condition and its two resulting error estimates are also established for the first time to our knowledge, for the robust block-sparse signal recovery using a regularized mixed $\ell_{2}/\ell_{1}$-norm minimization model, and these results well complement the existing theoretical investigation on this model which focuses on the non-uniform recovery conditions and/or the robust signal recovery in presence of the random noise.
\end{abstract}


%
\IEEEpeerreviewmaketitle

\section{Introduction}\label{Section1}
\IEEEPARstart{T}{HE} last decade has seen the burgeoning development of Compressed Sensing (CS) \cite{CandesRT-RUP-TIT-2006,Donoho-CS-TIT-2006} and its wide-spread applications in many fields. At the core of CS is the problem of efficiently recovering a sparse signal from a relatively small number of linear measurements. Mathematically, for any given signal $\bm{x}\in\mathbb{R}^{n}$, we say that it is sparse if most of its entries are zero. More specifically, if it has at most $k$ non-zero entries, i.e., $\|\bm{x}\|_{0}\triangleq |\text{supp}(\bm{x})|\leq k$, we call it a $k$-sparse signal. In standard CS, one usually observes the linear measurements of the sparse signal $\bm{x}$ via $\bm{b}=A\bm{x}$, where $A\in\mathbb{R}^{m\times n}$ ($m\ll n$) is a given measurement matrix. To recover such a sparse signal, a natural idea is to search the sparsest solution among all the possible solutions. This directly leads to the following $\ell_{0}$-norm minimization problem
\begin{align}\label{Constrained-L0}
\min_{\bm{x}\in\mathbb{R}^{n}}\|\bm{x}\|_{0},~~s.t.~~\bm{b}=A\bm{x}.
\end{align}
Unfortunately, this problem is NP-hard in general, and hence it is computationally infeasible. Instead, some algorithms which aims to pursue the suboptimal solutions of \eqref{Constrained-L0} were proposed, see, e.g., \cite{TroppG-OMP-TIT-2007,BlumensathD-IHT-ACHA-2009,Mohimani-etal-SLO-TSP-2009} and their variants \cite{NeedellV-ROMP-IEEESTSP-2010,Donoho-etal-StOMP-TIT-2012,KarahanogluErdogan-AOMP-DSP-2012,Foucart-HTP-SIAM-JNA-2011}. Importantly, many of these algorithms have been proved to perform well under certain conditions.

Besides the above algorithm strategies, there also exist many other efficient approaches \cite{Candes-Tao-L1-TIT-2005,Chartrand-Lp-SPL-2007,Candes-etal-WL1-JFAA-2008,Yin-etal-L12-SIAM-2015,LuLi-PSM-MOR-2018,Wang-IEEE-SPL-2017} which can circumvent the NP-hardness of \eqref{Constrained-L0}, and a popular one is the constrained $\ell_{1}$-norm minimization method which solves
\begin{align}\label{Noiseless-Constrained-L1}
\min_{\bm{x}\in\mathbb{R}^{n}}\|\bm{x}\|_{1},~~s.t.~~\bm{b}=A\bm{x},
\end{align}
where $\|\cdot\|_{1}$ is the $\ell_{1}$-norm of vector. Problem \eqref{Noiseless-Constrained-L1} is convex and therefore can be well addressed by many convex optimization softwares. To theoretically investigate the equivalence between \eqref{Constrained-L0} and \eqref{Noiseless-Constrained-L1}, one often adopts the \emph{Restricted Isometry Constant} (RIC) of matrix with $k$ order, denoted by $\delta_{k}$, which is defined to be the smallest value of $\delta\in(0, 1)$ such taht
\begin{align*}
(1-\delta)\|\bm{x}\|_{2}^{2}\leq\|A\bm{x}\|_{2}^{2}\leq(1+\delta)\|\bm{x}\|_{2}^{2}
\end{align*}
for every $k$-sparse vector $\bm{x}$. This notation was first proposed by Cand\`{e}s and Tao in \cite{Candes-Tao-L1-TIT-2005}, where they have shown that \eqref{Noiseless-Constrained-L1} is equivalent to \eqref{Constrained-L0} in noiselessly recovering any $k$-sparse signals when $\delta_{k}+\delta_{2k}+\delta_{3k}<1$. Subsequently many researchers were committed to improving this condition, see, e.g., \cite{Candes-RIP-CRM-2008,FoucartLai-RIC-ACHA-2009,Cai-etal-RIC-TSP-2010,MoLi-RIC-ACHA-2011,CaiZhang-RIC-ACHA-2013,TT-Cai-polytope-TIT-2014,R-Zhang-A-TIT-2018,Zhang-Li-Lp-ACHA-InPress}. In more general application scenarios of CS, one often wishes to recover the original signal $\bm{x}$ (may not be exactly sparse) from the noisy observation $\bm{b}$ with
\begin{align}\label{b=Ax+z}
\bm{b}=A\bm{x}+\bm{z},
\end{align}
where $\bm{z}\in\mathbb{R}^{n}$ is the unknown measurement noise/error, which directly leads to the following optimization problem
\begin{align}\label{Constrained-L1}
\min_{\bm{x}\in\mathbb{R}^{n}}\|\bm{x}\|_{1},~~s.t.~~\|\bm{b}-A\bm{x}\|_{2}\leq\epsilon,
\end{align}
where $\epsilon\geq0$. Obviously, \eqref{Constrained-L1} will reduces to \eqref{Noiseless-Constrained-L1} if one takes $\epsilon=0$. It should also be noted that the above-mentioned exact recovery conditions are still available to guarantee the robust recovery of signals from \eqref{Constrained-L1} in the presence of noise.

On the other hand, when used in many practical applications, particularly the applications where the input data are in large scale, the constrained problem \eqref{Constrained-L1} is not always convenient to solve. Instead, one often solves its unconstrained counterpart, i.e., the ($\ell_{2}$-norm) regularized $\ell_{1}$-norm minimization problem
\begin{align}\label{Unconstrained-L1-L2Noise}
\min_{\bm{x}\in\mathbb{R}^{n}}\|\bm{x}\|_{1}+\frac{1}{2\lambda}\|\bm{b}-A\bm{x}\|_{2}^{2},
\end{align}
where $\lambda$ is the nonnegative tradeoff parameter. This problem is also known as the Lasso estimator \cite{tibshirani1996regression} or the Basis Pursuit DeNoising (BPDN) \cite{Chen-BPDN-SIAM-2001} , and it can also be solved efficiently by many algorithms, see, e.g., \cite{Kim-etal-IPM-2007,GP-2007,FISTA-2009,ADMM-2011,ASBCD-2016}. Recently, the relation between \eqref{Constrained-L1} and \eqref{Unconstrained-L1-L2Noise} was carefully investigated by Zhang, Yuan and Yin \cite{Zhang-etal-Equality-2016} in the context of the \emph{non-uniform recovery} \cite{Foucart-Rauhut-CSBook-2013}, i.e., the recovery of some specific sparse signals, for example, the sparse signals limited in a specific support. However, when it comes to the \emph{uniform recovery} \cite{Zhang-etal-Equality-2016,Foucart-Rauhut-CSBook-2013}, i.e., recovering all the (general) sparse signals \footnote{In these signals, both the number and support of their non-zero entries are not known in advance.}, it is generally believed that \eqref{Constrained-L1} and \eqref{Unconstrained-L1-L2Noise} are not exact equivalents. As early as 2008, Zhu \cite{Zhu-1stRIC-Unconstrained-2008} has derived the RIC-based theoretical guarantee for \eqref{Unconstrained-L1-L2Noise}, which states that one can robustly recover any $k$-sparse signal $\bm{x}$ through \eqref{b=Ax+z} with $\|\bm{z}\|_{2}=\epsilon$ by using \eqref{Unconstrained-L1-L2Noise} under certain $\lambda$, if $A$ obeys $\delta_{4k}+2\delta_{5k}<1$. However, this work was relatively rarely noted by the researchers. In 2009, Bickel, Ritov and Tsybakov \cite{bickel2009simultaneous} established a RIC-like guarantee for \eqref{Unconstrained-L1-L2Noise} in the presence of the random noise. Recently, some new RIC conditions were obtained to ensure the robust recovery of some unconstrained $\ell_{1}$-analysis approaches under the Dantzig Selector (DS) type noise/error (i.e., $\|A^{T}\bm{z}\|_{\infty}\leq\epsilon$), see \cite{Lin-Li-ACHA-2014,Shen-etal-ACHA-2015,GeHM-UnconstrainedCS-Submitted-2018} for details.

Some random matrices represented by sub-Gaussian random matrices are proved to have a small RIC with overwhelmingly high probability \cite{Foucart-Rauhut-CSBook-2013}. However, when used in practical scenarios they often suffer from storage and computation limitations. Moreover, it is also NP-hard in general to find the RIC of any given matrix. To overcome these difficulties, some researchers proposed to reuse the \emph{mutual coherence}, another powerful tool introduced by Mallat and Zhang \cite{Mallat-Zhang-MP-Coherence-1993} in their initial research on matching pursuit. In this paper, denoting by $\bm{a}_{i}\in\mathbb{R}^{m}$ the $i$th column of matrix $A$, we shall define the mutual coherence of matrix $A\triangleq[\bm{a}_{1},\bm{a}_{2},\cdots, \bm{a}_{n}]$ as
\begin{align*}
\mu=\max_{1\leq i< j\leq n}|\langle\bm{a}_{i},\bm{a}_{j}\rangle|,
\end{align*}
where we assume that $\bm{a}_{i}$ obeys $\|\bm{a}_{i}\|_{2}=1$ for $i=1,2,\cdots, n$. Many deterministic measurement matrices in fact are designed according to the mutual coherence. There are many coherence-based theoretical results for \eqref{Constrained-L1}, see, e.g., \cite{Donoho-Huo-Coherence-TIT-2001,GribonvalNielsen-Coherence-TIT-2003,Fuchs-Coherence-TIT-2004,Donoho-etal-NoiseCoherence-TIT-2006,Cai-etal-coherence-TIT-2009,Cai-etal-SharpCoherence-TIT-2010}. In particular in \cite{Cai-etal-SharpCoherence-TIT-2010}, Cai, Wang and Xu have shown that any signals (not necessary to be $k$-sparse) can be robustly recovered using \eqref{Constrained-L1}, if $A$ satisfies $\mu<1/(2k-1)$, and this condition is also sharp for the noiseless recovery of any $k$-sparse signals through \eqref{Noiseless-Constrained-L1}.

As far as we know, the first (mutual) coherence-based result on \eqref{Unconstrained-L1-L2Noise} was given by Fuchs \cite{Fuchs-Coherence-TIT-2004} in 2004 under the non-uniform recovery setting, which states that any fixed signal $\bm{x}$ with $k$ non-zero entries (i.e., $\|\bm{x}\|_{0}=k$) can be uniquely recovered from $\bm{b}=A\bm{x}$ using \eqref{Unconstrained-L1-L2Noise}, if $A$ satisfies $\mu<1/(2k-1)$ and $\lambda$ in \eqref{Unconstrained-L1-L2Noise} has been taken small enough. Later, Fuchs \cite{fuchs2005recovery} further investigated \eqref{Unconstrained-L1-L2Noise} for the noisy signal recovery, and shown that if $\bm{b}$ is observed through \eqref{b=Ax+z} with $\|\bm{x}\|_{0}=k$ and $\|\bm{z}\|_{2}\leq\epsilon$,  $A$ satisfies $\mu\leq c/k$ for certain $c\leq1/2$, and $\bm{x}^{\sharp}$ is assumed to be the optimal solution of \eqref{Unconstrained-L1-L2Noise} under certain $\lambda$ (related to $\mu$, $c$ and $\epsilon$), then the support of $\bm{x}^{\sharp}$ will be either identical to, or contained in, that of $\bm{x}$. Moreover, Fuchs also shown that if similar constraints are imposed on $\bm{x}$, $\bm{z}$ and $\lambda$, $\bm{x}^{\sharp}$ and $\bm{x}$ will have their non-zero entries at the same support and with the same signs. Subsequently, Tropp \cite{tropp2006just} further extended the results in \cite{fuchs2005recovery} to more general case. In 2010, Ben-Haim, Eldar and Elad \cite{ben2010coherence} revisited \eqref{Unconstrained-L1-L2Noise} under the random noise, and their obtained coherence results have been proved to be better than those induced in \cite{bickel2009simultaneous}. Note that the above coherence results apply only to deal with the signals whose sparsity is known in advance. Recently, using the \emph{cumulative coherence} \cite{tropp2004greed} tool, Li and Chen \cite{li2019signal} established a new uniform recovery condition for \eqref{Unconstrained-L1-L2Noise} to deal with the signal recovery in the presence of noise. Their results shown that if $A$ obeys $\mu\leq1/(\sqrt{3}(5k-2))$, one can robustly recover any signals corrupted with the DS type noise. However, the noise they considered is based on the DS type noise rather than the often used $\ell_{2}$-norm bounded noise, and the recovery condition they obtained still has much room to improve.

In this paper, equipped with the powerful mutual coherence tool, we investigate the performance guarantees of \eqref{Unconstrained-L1-L2Noise} and its some variants. In summary, the contributions of this paper are listed as follows:
\begin{itemize}
  \item We establish a tight uniform recovery condition and two relatively tight error estimates for \eqref{Unconstrained-L1-L2Noise}, which are sufficient to guarantee the robust recovery of signals corrupted with the $\ell_{2}$-norm bounded noise.
  \item We extend the obtained theoretical results to guarantee the robust recovery of the signals corrupted with the DS type noise and also the structured block-sparse signal recovery in the presence of the bounded noise. To the best of our knowledge, these extended results are established for the first time under the uniform recovery setting.
\end{itemize}

The remainder of this paper is organized as follows. Section \ref{Notations-Preliminaries} introduces some notations and preliminaries. In Section \ref{Main-results}, we present the main results. Section \ref{Extensions} shows two extensions. Finally, conclusion and future work are given in Section \ref{Conclusions-Futurework}.

\section{Notations and Preliminaries}
\label{Notations-Preliminaries}

\subsection{Notations}
Throughout this paper, we denote $[r]\triangleq\{1,2,\cdots, r\}$ for any given integer $r$, and $E^{c}=[n]\setminus E$ for any given index set $E\subset[n]$. We also denote $\bm{h}_{E}$ as a vector whose entries $(\bm{h}_{E})_{i}=\bm{h}_{i}$ for $i\in E$ and 0 otherwise, and $\|\cdot\|_{a}^{b}=(\|\cdot\|_{a})^{b}$ where $\|\cdot\|_{a}$ represents certain norm or quasi-norm. For any signal $\bm{x}\in\mathbb{R}^{n}$, we denote its best $k$-term approximate as $\bm{x}_{[k]}$, which is defined as
\begin{align*}
\bm{x}_{[k]}=\mathop{\arg\min}_{\|\bm{y}\|_{0}\leq k}\|\bm{y}-\bm{x}\|_{2}.
\end{align*}
Besides, for the simplicity of symbol expression we introduce the following two functions
\begin{align*}
f_{a}(x)=ax^{2}+3\sqrt{a}x+3,~~g_{a}(x)=2ax^{2}+4\sqrt{a}x+1.
\end{align*}

\subsection{Three key lemmas}
The proof of our main results heavily relies on the following three lemmas. We start with introducing the first lemma, which provides a RIC-like coherence result for any given matrix.

\begin{lemma}[\cite{tropp2006just,Donoho-etal-NoiseCoherence-TIT-2006}]\label{RIC-like-coherence-Lemma}
Assume that the matrix $A\in\mathbb{R}^{m\times n}$ is standardized to have unit $\ell_{2}$-norm. Then it holds that
\begin{align}\label{quasi-RIP-Coherence}
(1-(k-1)\mu)\|\bm{y}\|_{2}^{2}\leq\|A\bm{y}\|_{2}^{2}\leq(1+(k-1)\mu)\|\bm{y}\|_{2}^{2}
\end{align}
for all $k$-sparse signals $\bm{y}\in\mathbb{R}^{n}$.
\end{lemma}

\begin{lemma}\label{NP-Lemma1}
If $\bm{b}$ is observed via \eqref{b=Ax+z} with $\|\bm{z}\|_{2}\leq\epsilon$, then for any subset $E\subset[n]$ with $|E|=k$, the optimal solution
$\bm{x}^{\sharp}$ of \eqref{Unconstrained-L1-L2Noise} satisfies
\begin{align}\label{NP-Lemma1-results1}
\|A\bm{h}\|_{2}^{2}-2\epsilon\|A\bm{h}\|_{2}\leq&2\lambda(\|\bm{h}_{E}\|_{1}-\|\bm{h}_{E^{c}}\|_{1}+2\|\bm{x}_{E^{c}}\|_{1})
\end{align}
and
\begin{align}\label{NP-Lemma1-results2}
\|\bm{h}_{E^{c}}\|_{1}\leq\|\bm{h}_{E}\|_{1}+2\|\bm{x}_{E^{c}}\|_{1}+\frac{\epsilon}{\lambda}\|A\bm{h}\|_{2},
\end{align}
where $\bm{h}=\bm{x}^{\sharp}-\bm{x}$.
\end{lemma}

\begin{IEEEproof}[Proof of Lemma \ref{NP-Lemma1}]
Since $\bm{x}^{\sharp}$ is the optimal solution of \eqref{Unconstrained-L1-L2Noise}, we have
\begin{align*}
\|\bm{x}^{\sharp}\|_{1}+\frac{1}{2\lambda}\|\bm{b}-A\bm{x}^{\sharp}\|_{2}^{2}\leq\|\bm{x}\|_{1}+\frac{1}{2\lambda}\|\bm{b}-A\bm{x}\|_{2}^{2},
\end{align*}
which is equivalent to
\begin{align}\label{Unconstrained-L1-L2Noise-Estimate-LR}
\|A\bm{h}\|_{2}^{2}-2\langle \bm{z}, A\bm{h}\rangle\leq 2\lambda(\|\bm{x}\|_{1}-\|\bm{x}^{\sharp}\|_{1}).
\end{align}
As to the Left-Hand Side (LHS) of \eqref{Unconstrained-L1-L2Noise-Estimate-LR}, we have
\begin{align}\label{Unconstrained-L1-L2Noise-Estimate-L}
\text{LHS}\geq\|A\bm{h}\|_{2}^{2}-2\epsilon\|A\bm{h}\|_{2}.
\end{align}
As to the Right-Hand Side (RHS) of \eqref{Unconstrained-L1-L2Noise-Estimate-LR}, we know
\begin{align}\label{Unconstrained-L1-L2Noise-Estimate-R}
\text{RHS}=&\sum_{i=1}^{n}|\bm{x}_{i}+\bm{h}_{i}|-\|\bm{x}\|_{1}\nonumber\\      %
\geq&\sum_{i\in E}|\bm{x}_{i}-\bm{h}_{i}|+\sum_{i\in E^{c}}|\bm{h}_{i}-\bm{x}_{i}|-\|\bm{x}\|_{1}\nonumber\\
\geq&\sum_{i\in E}(|\bm{x}_{i}|-|\bm{h}_{i}|)+\sum_{i\in E^{c}}(|\bm{h}_{i}|-|\bm{x}_{i}|)-(\|\bm{x}_{E}\|_{1}+\|\bm{x}_{E^{c}}\|_{1})\nonumber\\
=&-\|\bm{h}_{E}\|_{1}+\|\bm{h}_{E^{c}}\|_{1}-2\|\bm{x}_{E^{c}}\|_{1}.
\end{align}
Therefore combing \eqref{Unconstrained-L1-L2Noise-Estimate-LR}, \eqref{Unconstrained-L1-L2Noise-Estimate-L} and \eqref{Unconstrained-L1-L2Noise-Estimate-R} directly leads to \eqref{NP-Lemma1-results1}, and \eqref{NP-Lemma1-results2} follows trivially from \eqref{NP-Lemma1-results1}.
\end{IEEEproof}

\begin{lemma}\label{NP-Lemma2}
If the matrix $A\in\mathbb{R}^{m\times n}$ is standardized to have unit $\ell_{2}$-norm, and obeys
\begin{align}\label{Condtion-NP-Lemma2}
\mu<\frac{1}{k-1}
\end{align}
for certain integer $k\geq2$, then for any vector $\bm{h}\in\mathbb{R}^{n}$ and any subset $E\subset[n]$ with $|E|=k$, it holds that
\begin{align}\label{NP-Lemma2-results}
\|\bm{h}_{E}\|_{2}\leq \alpha_{1}\|A\bm{h}\|_{2}+\alpha_{2}\|\bm{h}_{E^{c}}\|_{1},
\end{align}
where
\begin{align*}
\alpha_{1}\triangleq\frac{\sqrt{1+(k-1)\mu}}{1-(k-1)\mu} \text{~~and~~} \alpha_{2}\triangleq \frac{\sqrt{k}\mu}{1-(k-1)\mu}.
\end{align*}
\end{lemma}

\begin{remark}
It is easy to know from Lemma \ref{NP-Lemma2} that both $\alpha_{1}$ and $\alpha_{2}$ are two monotone increasing functions on variable $\mu$. Therefore if one restricts $\mu<1/(2k-1)$, it will be clear that
\begin{align}\label{Simple-Calculation}
1<\alpha_{1}<\sqrt{6}, \sqrt{k}\alpha_{2}<1 \text{~and~}\frac{1}{1-\sqrt{k}\alpha_{2}}<\frac{1}{1-(2k-1)\mu}.
\end{align}
\end{remark}

\begin{IEEEproof}[Proof of Lemma \ref{NP-Lemma2}]
The proof is simple. We start with estimating the lower and upper bounds of
\begin{align*}
\rho\triangleq|\langle A\bm{h},A\bm{h}_{E}\rangle|.
\end{align*}
First, using Lemma \ref{RIC-like-coherence-Lemma}, we know
\begin{align}\label{NP-Lemma2s-Lower}
\rho\geq&|\langle A\bm{h}_{E},A\bm{h}_{E}\rangle|-|\langle A\bm{h}_{E^{c}},A\bm{h}_{E}\rangle|\nonumber\\
\geq&(1-(k-1)\mu)\|\bm{h}_{E}\|_{2}^{2}-|\sum_{i\in E}\sum_{j\in E^{c}}\langle\bm{a}_{i},\bm{a}_{j}\rangle\bm{h}_{i}\bm{h}_{j}|\nonumber\\
\geq&(1-(k-1)\mu)\|\bm{h}_{E}\|_{2}^{2}-\mu\|\bm{h}_{E}\|_{1}\|\bm{h}_{E^{c}}\|_{1}\nonumber\\
\geq&(1-(k-1)\mu)\|\bm{h}_{E}\|_{2}^{2}-\sqrt{k}\mu\|\bm{h}_{E}\|_{2}\|\bm{h}_{E^{c}}\|_{1}.
\end{align}
Next we estimate the upper bound of $\rho$. It follows from \ref{RIC-like-coherence-Lemma} that
\begin{align}\label{NP-Lemma2s-Upper}
\rho\leq&\|A\bm{h}\|_{2}\|A\bm{h}_{E}\|_{2}\leq\sqrt{1+(k-1)\mu}\|\bm{h}_{E}\|_{2}\|A\bm{h}\|_{2}.
\end{align}
Now, combing \eqref{NP-Lemma2s-Lower}, \eqref{NP-Lemma2s-Upper} and the condition \eqref{Condtion-NP-Lemma2} directly leads to the desired inequality \eqref{NP-Lemma2-results}.
\end{IEEEproof}

\section{Main Results}
\label{Main-results}

With preparations above, we now present our main results.
\begin{theorem}\label{Theorem-L2Noise}
For any $\bm{b}$ observed via \eqref{b=Ax+z} with $\|\bm{z}\|_{2}\leq\epsilon$, if the measurement matrix $A$, which is standardized to have unit $\ell_{2}$-norm, satisfies
\begin{align}\label{sharp-mutual-coherence-L1}
\mu<\frac{1}{2k-1}
\end{align}
for certain integer $k\geq2$, then we have
\begin{align}\label{Theorem-L2Noise-results1}
\|A(\bm{x}^{\sharp}-\bm{x})\|_{2}&\leq C_{1}(\alpha_{1})\|\bm{x}-\bm{x}_{[k]}\|_{1} + C_{2}(\alpha_{1}), \\
\|\bm{x}^{\sharp}-\bm{x}\|_{2}&\leq C_{3}(\alpha_{1},\alpha_{2})\|\bm{x}-\bm{x}_{[k]}\|_{1}+ C_{4}(\alpha_{1},\alpha_{2}),\nonumber
\end{align}
where $\bm{x}^{\sharp}$ is the optimal solution of \eqref{Unconstrained-L1-L2Noise}, and
\begin{align*}
C_{1}(\alpha_{1})&\triangleq\frac{2\lambda}{\sqrt{k}\alpha_{1}\lambda+\epsilon},\\
C_{2}(\alpha_{1})&\triangleq2(\sqrt{k}\alpha_{1}\lambda+\epsilon),\\
C_{3}(\alpha_{1},\alpha_{2})&\triangleq\frac{2\sqrt{k}\alpha_{1}f_{k}(\alpha_{2})\lambda+2g_{k}(\alpha_{2})\epsilon}{\sqrt{k}(1-\sqrt{k}\alpha_{2})(\sqrt{k}\alpha_{1}\lambda+\epsilon)},\\
C_{4}(\alpha_{1},\alpha_{2})&\triangleq\frac{\sqrt{k}\alpha_{1}(5+2\sqrt{k}\alpha_{2})\lambda+g_{k}(\alpha_{2})\epsilon}{\sqrt{k}(1-\sqrt{k}\alpha_{2})\lambda(\sqrt{k}\alpha_{1}\lambda+\epsilon)^{-1}}.
\end{align*}
\end{theorem}
\begin{remark}
Theorem \ref{Theorem-L2Noise} shows that one can robustly recover any signals (may not be $k$-sparse) corrupted with the $\ell_{2}$-norm bounded noise, if the measurement matrix $A$ satisfies \eqref{sharp-mutual-coherence-L1}. To the best of our knowledge, we first extend this sharp uniform recovery condition \footnote{The sharp condition/bound here and throughout will refer to $\mu<1/(2k-1)$. It has been shown in \cite{Cai-etal-SharpCoherence-TIT-2010} that, any $k$-sparse signal $\bm{x}$, without exception, can be exactly recovered from \eqref{Noiseless-Constrained-L1} under this sharp condition, and there exists a matrix $A$ with $m<n$ obeying $u=1/(2k-1)$, and two nonzero $k$-sparse vectors $\bm{\widetilde{x}}$ and $\bm{\widehat{x}}$ with disjoint supports such that $A\bm{\widetilde{x}}\neq A\bm{\widehat{x}}$.} derived by Cai, Wang and Xu in \cite{Cai-etal-SharpCoherence-TIT-2010} for the constrained problem \eqref{Noiseless-Constrained-L1} to its unconstrained counterpart, i.e., the unconstrained problem \eqref{Unconstrained-L1-L2Noise}. Similar to \cite{Shen-etal-ACHA-2015,GeHM-UnconstrainedCS-Submitted-2018}, if we associate $\epsilon$ with $\lambda$, e.g., setting $\epsilon=\lambda$, then we get a special case of Theorem \ref{Theorem-L2Noise}, and one can find this result in Corollary \ref{corollary1}. In Remark \ref{remark-lambda-epsilon}, we will analyze the tightness of these two error estimates under the setting of $\epsilon=\lambda$. Besides, it is also very easy to induce some other special cases of Theorem \ref{Theorem-L2Noise} to cope with several different sparse recovery tasks. For examples, one can consider the robust recovery of any exactly sparse signals, i.e., setting the original signals $\bm{x}$ to be exactly $k$-sparse. The detailed analysis of these cases will become very similar to that of Corollary \ref{corollary1}.
\end{remark}
\begin{corollary}\label{corollary1}
Assume that $\bm{b}$ is observed via \eqref{b=Ax+z} with $\|\bm{z}\|_{2}\leq\lambda$. If the measurement matrix $A$ is standardized to have unit $\ell_{2}$-norm, and also satisfies \eqref{sharp-mutual-coherence-L1} for certain integer $k\geq2$, then we have
\begin{align*}
\|A(\bm{x}^{\sharp}-\bm{x})\|_{2}&\leq \widehat{C}_{1}\|\bm{x}-\bm{x}_{[k]}\|_{1} + \widehat{C}_{2}\lambda,\\
\|\bm{x}^{\sharp}-\bm{x}\|_{2}&\leq \widehat{C}_{3}\|\bm{x}-\bm{x}_{[k]}\|_{1}+ \widehat{C}_{4}\lambda,
\end{align*}
where $\bm{x}^{\sharp}$ is the optimal solution of \eqref{Unconstrained-L1-L2Noise}, and
\begin{align*}
\widehat{C}_{1}&=\frac{2}{\sqrt{k}\alpha_{1}+1},\\
\widehat{C}_{2}&=2(\sqrt{k}\alpha_{1}+1),\\
\widehat{C}_{3}&=\frac{2\sqrt{k}\alpha_{1}f_{k}(\alpha_{2})+2g_{k}(\alpha_{2})}{\sqrt{k}(1-\sqrt{k}\alpha_{2})(\sqrt{k}\alpha_{1}+1)},\\
\widehat{C}_{4}&=\frac{\sqrt{k}\alpha_{1}(5+2\sqrt{k}\alpha_{2})+g_{k}(\alpha_{2})}{\sqrt{k}(1-\sqrt{k}\alpha_{2})(\sqrt{k}\alpha_{1}+1)^{-1}}.
\end{align*}
\end{corollary}

\begin{remark}\label{remark-lambda-epsilon}
Due to the existence of $\alpha_{1}$ and $\alpha_{2}$, the coefficients $\widehat{C}_{1}, \cdots, \widehat{C}_{4}$ are not convenient to be analyzed. Fortunately, based on the previous estimates for $\alpha_{1}$ and $\alpha_{2}$, i.e., \eqref{Simple-Calculation}, we can give a rough but simple estimate for each $\widehat{C}_{i}$. Specifically,
\begin{align*}
\widehat{C}_{1}<&\frac{2}{\sqrt{k}+1}<\frac{2}{\sqrt{k}},\\
\widehat{C}_{2}<&2(\sqrt{6k}+1)\leq2(\sqrt{6}+1)\sqrt{k},\\
\widehat{C}_{3}<&\frac{14(\sqrt{6k}+1)}{\sqrt{k}(\sqrt{k}+1)(1-\sqrt{k}\alpha_{2})}<\frac{14\sqrt{6}}{\sqrt{k}\left(1-(2k-1)\mu\right)},\\
\widehat{C}_{4}<&\frac{7(\sqrt{6k}+1)^{2}}{\sqrt{k}(1-\sqrt{k}\alpha_{2})}<\frac{7(\sqrt{6}+1)^{2}\sqrt{k}}{1-(2k-1)\mu}.
\end{align*}
These upper bound estimates of coefficients make our recovery error, denoted by $\text{RE}$, have the form of
\begin{align}\label{public-form}
\text{RE}\leq \overline{C}_{1}\frac{\|\bm{x}-\bm{x}_{[k]}\|_{1}}{\sqrt{k}} + \overline{C}_{2}\sqrt{k}\lambda,
\end{align}
where $\text{RE}$ stands for $\|A(\bm{x}^{\sharp}-\bm{x})\|_{2}$ or $\|\bm{x}^{\sharp}-\bm{x}\|_{2}$, and $\overline{C}_{1}$ and $\overline{C}_{2}$ depend only on the value of $1-(2k-1)\mu$, which characterizes the gap between the coherence of the selected measurement matrix $A$ and its sharp bound. This result also coincides with the ones obtained in \cite{Lin-Li-ACHA-2014,Shen-etal-ACHA-2015,GeHM-UnconstrainedCS-Submitted-2018,li2019signal} for \eqref{Unconstrained-L1-L2Noise} in form. However, one should note that the authors in these literature focus on sparse recovery corrupted with the DS type noise, which is totally different from ours. Despite this, our obtained upper bound estimates to some degree are still better than theirs since a much tighter (or sharp) recovery condition is used. What's more, some  coefficients in these estimates can be further improved if one optimizes some inequalities used to prove Theorem \ref{Theorem-L2Noise}.
\end{remark}

\begin{IEEEproof}[Proof of Theorem \ref{Theorem-L2Noise}]
We first assume that $\mu<1/(k-1)$ for certain integer $k\geq2$ and denote
\begin{align*}
E=\text{supp}(\bm{x}_{[k]}) \text{~~and~~} \bm{h}=\bm{x}^{\sharp}-\bm{x},
\end{align*}
then using Lemma \ref{NP-Lemma1}, Lemma \ref{NP-Lemma2}, we have
\begin{align}\label{Theorem-L2Noise-Proof-1}
\|A\bm{h}\|_{2}^{2}-2\epsilon\|A\bm{h}\|_{2}\leq&2\lambda(\|\bm{h}_{E}\|_{1}-\|\bm{h}_{E^{c}}\|_{1}+2\|\bm{x}_{E^{c}}\|_{1})\nonumber\\
\leq&2\sqrt{k}\lambda(\alpha_{1}\|A\bm{h}\|_{2}+\alpha_{2}\|\bm{h}_{E^{c}}\|_{1})-2\lambda\|\bm{h}_{E^{c}}\|_{1}+4\lambda\|\bm{x}_{E^{c}}\|_{1}\nonumber\\
=&2\sqrt{k}\alpha_{1}\lambda\|A\bm{h}\|_{2}+4\lambda\|\bm{x}_{E^{c}}\|_{1}-2(1-\sqrt{k}\alpha_{2})\lambda\|\bm{h}_{E^{c}}\|_{1},
\end{align}
where we used $\|\bm{h}_{E}\|_{1}\leq\sqrt{k}\|\bm{h}_{E}\|_{2}$. We can known from the condition \eqref{sharp-mutual-coherence-L1} that
\begin{align*}
1-\sqrt{k}\alpha_{2}&=1-\frac{k\mu}{1-(k-1)\mu}\\
&>1-\frac{k/(2k-1)}{1-(k-1)/(2k-1)}=0.
\end{align*}
Therefore we can further write \eqref{Theorem-L2Noise-Proof-1} as
\begin{align*}
\|A\bm{h}\|_{2}^{2}-2(\epsilon+\sqrt{k}\alpha_{1}\lambda)\|A\bm{h}\|_{2}-4\lambda\|\bm{x}_{E^{c}}\|_{1}\leq0.
\end{align*}
which implies that
\begin{align*}
\|A\bm{h}\|_{2}\leq&(\sqrt{k}\alpha_{1}\lambda+\epsilon)+\sqrt{(\sqrt{k}\alpha_{1}\lambda+\epsilon)^{2}+ 4\lambda\|\bm{x}_{E^{c}}\|_{1}}\\
\leq&(\sqrt{k}\alpha_{1}\lambda+\epsilon)+\sqrt{\left(\sqrt{k}\alpha_{1}\lambda+\epsilon+\frac{2\lambda\|\bm{x}_{E^{c}}\|_{1}}{\sqrt{k}\alpha_{1}\lambda+\epsilon}\right)^{2}}\\
=&\frac{2\lambda}{\sqrt{k}\alpha_{1}\lambda+\epsilon}\|\bm{x}_{E^{c}}\|_{1}+2(\sqrt{k}\alpha_{1}\lambda+\epsilon).
\end{align*}
This completes \eqref{Theorem-L2Noise-results1}. Based on \eqref{Theorem-L2Noise-results1}, \eqref{NP-Lemma1-results2} and \eqref{NP-Lemma2-results}, we have
\begin{align}\label{NP-Lemma1-result2-new}
\|\bm{h}_{E}\|_{2}\leq& \alpha_{1}\|A\bm{h}\|_{2}+\alpha_{2}\bigg(\|\bm{h}_{E}\|_{1}+2\|\bm{x}_{E^{c}}\|_{1}+\frac{\epsilon}{\lambda}\|A\bm{h}\|_{2}\bigg)\nonumber\\
=&\frac{\alpha_{1}\lambda+\alpha_{2}\epsilon}{\lambda}\|A\bm{h}\|_{2}+\sqrt{k}\alpha_{2}\|\bm{h}_{E}\|_{2}+2\alpha_{2}\|\bm{x}_{E^{c}}\|_{1}\nonumber\\
\leq&\frac{\alpha_{1}\lambda+\alpha_{2}\epsilon}{(1-\sqrt{k}\alpha_{2})\lambda}\|A\bm{h}\|_{2}+\frac{2\alpha_{2}}{1-\sqrt{k}\alpha_{2}}\|\bm{x}_{E^{c}}\|_{1}\nonumber\\
\leq&\frac{\alpha_{1}\lambda+\alpha_{2}\epsilon}{(1-\sqrt{k}\alpha_{2})\lambda}\bigg(\frac{2\lambda}{\sqrt{k}\alpha_{1}\lambda+\epsilon}\|\bm{x}_{E^{c}}\|_{1}+2\sqrt{k}\alpha_{1}\lambda
+2\epsilon\bigg)+\frac{2\alpha_{2}}{1-\sqrt{k}\alpha_{2}}\|\bm{x}_{E^{c}}\|_{1}\nonumber\\
=&\frac{2\alpha_{1}(1+\sqrt{k}\alpha_{2})\lambda+4\alpha_{2}\epsilon}{(1-\sqrt{k}\alpha_{2})(\sqrt{k}\alpha_{1}\lambda+\epsilon)}\|\bm{x}_{E^{c}}\|_{1}
+\frac{2(\alpha_{1}\lambda+\alpha_{2}\epsilon)(\sqrt{k}\alpha_{1}\lambda+\epsilon)}{(1-\sqrt{k}\alpha_{2})\lambda}.
\end{align}
Besides, using \eqref{NP-Lemma1-result2-new} together with \eqref{NP-Lemma1-results2} and \eqref{Theorem-L2Noise-results1} again, we can estimate $\|\bm{h}_{E^{c}}\|_{1}$ as
\begin{align}\label{NP-Lemma1-results2-new}
\|\bm{h}_{E^{c}}\|_{1}\leq\frac{4}{1-\sqrt{k}\alpha_{2}}\|\bm{x}_{E^{c}}\|_{1}
+\frac{2(\sqrt{k}\alpha_{1}\lambda+\epsilon)^{2}}{(1-\sqrt{k}\alpha_{2})\lambda}.
\end{align}
On the other hand, let $E_{1}$ be the index set of the $k$ largest entries of $\bm{h}_{E^{c}}$. Then we know from \cite{Shen-etal-ACHA-2015} that
\begin{align}\label{2-3-Shenyi}
\|\bm{h}_{E^{c}}\|_{2}\leq\|\bm{h}_{E_{1}}\|_{2}+\frac{\|\bm{h}_{E^{c}}\|_{1}}{2\sqrt{k}}.
\end{align}
Similarly, using Lemma \ref{NP-Lemma2} again on index $E_{1}$, we also have
\begin{align*}%
\|\bm{h}_{E_{1}}\|_{2}\leq \alpha_{1}\|A\bm{h}\|_{2}+\alpha_{2}\|\bm{h}_{E_{1}^{c}}\|_{1}.
\end{align*}
This, together with \eqref{Theorem-L2Noise-results1} and \eqref{NP-Lemma1-result2-new}, directly leads to
\begin{align}\label{NP-Lemma2-results-E1}
\|\bm{h}_{E_{1}}\|_{2}\leq& \alpha_{1}\left(\frac{2\lambda}{\sqrt{k}\alpha_{1}\lambda+\epsilon}\|\bm{x}_{E^{c}}\|_{1}+2\sqrt{k}\alpha_{1}\lambda+2\epsilon\right)+\alpha_{2}(\|\bm{h}_{E}\|_{1}+\|\bm{h}_{E^{c}}\|_{1})\nonumber\\
\leq&\frac{2\alpha_{1}\left(1+k(\alpha_{2})^{2}\right)\lambda+4\sqrt{k}(\alpha_{2})^{2}\epsilon}{(1-\sqrt{k}\alpha_{2})(\sqrt{k}\alpha_{1}\lambda+\epsilon)}\|\bm{x}_{E^{c}}\|_{1}
+\frac{2\left(\alpha_{1}\lambda+\sqrt{k}(\alpha_{2})^{2}\epsilon\right)}{(1-\sqrt{k}\alpha_{2})\lambda(\sqrt{k}\alpha_{1}\lambda+\epsilon)^{-1}}+\alpha_{2}\|\bm{h}_{E^{c}}\|_{1}.
\end{align}
Now, combining \eqref{NP-Lemma1-result2-new}, \eqref{NP-Lemma1-results2-new} and \eqref{NP-Lemma2-results-E1}, we can estimate $\|\bm{h}\|_{2}$ as follows:
\begin{align*}
\|\bm{h}\|_{2}\leq&\|\bm{h}_{E}\|_{2}+\|\bm{h}_{E^{c}}\|_{2}\\
\leq&\|\bm{h}_{E}\|_{2}+\|\bm{h}_{E_{1}}\|_{2}+\frac{\|\bm{h}_{E^{c}}\|_{1}}{2\sqrt{k}}\\
\leq&\frac{2\sqrt{k}\alpha_{1}f_{k}(\alpha_{2})\lambda+2g_{k}(\alpha_{2})\epsilon}{\sqrt{k}(1-\sqrt{k}\alpha_{2})(\sqrt{k}\alpha_{1}\lambda+\epsilon)}\|\bm{x}_{E^{c}}\|_{1}
+\frac{\sqrt{k}\alpha_{1}(5+2\sqrt{k}\alpha_{2})\lambda+2g_{k}(\alpha_{2})\epsilon}{\sqrt{k}(1-\sqrt{k}\alpha_{2})\lambda(\sqrt{k}\alpha_{1}\lambda+\epsilon)^{-1}},
\end{align*}
which completes the proof.
\end{IEEEproof}

\section{Extensions}
\label{Extensions}

In this section, two extensions of Theorem \ref{Theorem-L2Noise} are discussed. They include extending Theorem \ref{Theorem-L2Noise} to guarantee the robust recovery of signals from a DS regularized $\ell_{1}$-norm minimization model in the presence of the DS type noise, and that of the structured block-sparse signals from two regularized mixed $\ell_{2}/\ell_{1}$-norm minimization models in the presence of the bounded noise. We start with introducing the DS regularized $\ell_{1}$-norm minimization model for signal recovery in the presence of the DS type noise.

\subsection{Robust recovery via a DS regularized $\ell_{1}$-norm minimization}

The research on the DS type noise was initiated by Cand\`{e}s and Tao in \cite{candes2007dantzig}, which aims at recovering the signals corrupted with the DS type noise by solving the following constrained problem
\begin{align}\label{Constrained-L1-DS}
\min_{\bm{x}\in\mathbb{R}^{n}}\|\bm{x}\|_{1},~~s.t.~~\|A^{T}(\bm{b-A\bm{x}})\|_{\infty}\leq\epsilon.
\end{align}
Many remarkable results on this problem have been obtained over the past decade, see, e.g., \cite{bickel2009simultaneous,Cai-etal-RIC-TSP-2010,CaiZhang-RIC-ACHA-2013,TT-Cai-polytope-TIT-2014,R-Zhang-A-TIT-2018,Zhang-Li-Lp-ACHA-InPress,Cai-etal-SharpCoherence-TIT-2010}. Similar to the relation of \eqref{Constrained-L1} and \eqref{Unconstrained-L1-L2Noise}, a closely related problem to \eqref{Constrained-L1-DS} is the following DS regularized $\ell_{1}$-norm minimization problem:
\begin{align}\label{unconstrained-L1-DS}
\min_{\bm{x}\in\mathbb{R}^{n}}\|\bm{x}\|_{1}+\frac{1}{2\lambda}\|A^{T}(\bm{b-A\bm{x}})\|_{\infty}^{2}.
\end{align}
Inspired by Theorem \ref{Theorem-L2Noise}, we also establish a uniform recovery condition and two relatively tight error estimates for \eqref{unconstrained-L1-DS} to guarantee the robust signal recovery in the presence of such kind of noise, see Theorem \ref{Theorem-L2Noise-Extension-1} for details. This new theorem as well as Theorem \ref{Theorem-L2Noise}, to the best of our knowledge, first extends the sharp uniform recovery condition obtained in \cite{Cai-etal-SharpCoherence-TIT-2010} for \eqref{Noiseless-Constrained-L1} to its two unconstrained variants, i.e., \eqref{Unconstrained-L1-L2Noise} and \eqref{unconstrained-L1-DS}, to deal with the signals corrupted with the $\ell_{2}$-norm bounded noise and the DS type noise, respectively. In what follows, we present this theorem.

\begin{theorem}\label{Theorem-L2Noise-Extension-1}
For any $\bm{b}$ observed via \eqref{b=Ax+z} with $\|A^{T}\bm{z}\|_{2}\leq\epsilon$, if the measurement matrix $A$, whose columns are standardized to have unit $\ell_{2}$-norm, satisfies
\begin{align}\label{Theorem-L2Noise-condition1-Extension-2}
\mu<\frac{1}{2k-1}
\end{align}
for certain integer $k\geq2$, then we have
\begin{align*}
\|A^{T}A(\bm{x}^{\sharp}-\bm{x})\|_{\infty}&\leq C_{1}(\alpha_{1})\|\bm{x}-\bm{x}_{[k]}\|_{1} + C_{2}(\alpha_{1}), \\
\|\bm{x}^{\sharp}-\bm{x}\|_{2}&\leq C_{3}(\alpha_{1},\alpha_{2})\|\bm{x}-\bm{x}_{[k]}\|_{1}+ C_{4}(\alpha_{1},\alpha_{2}),\nonumber
\end{align*}
where $\bm{x}^{\sharp}$ here denotes the optimal solution of \eqref{unconstrained-L1-DS}.
\end{theorem}
\begin{remark}
In general, it is usually suggested to recover the signals corrupted with the DS type noise using the constrained problem \eqref{Constrained-L1-DS}, see, e.g., \cite{candes2007dantzig,TT-Cai-polytope-TIT-2014,R-Zhang-A-TIT-2018,Zhang-Li-Lp-ACHA-InPress}. Recently, some researchers proposed to deal with such kind of noise using the unconstrained problem \eqref{Unconstrained-L1-L2Noise}, see, e.g., \cite{Lin-Li-ACHA-2014,Shen-etal-ACHA-2015,GeHM-UnconstrainedCS-Submitted-2018,li2019signal}, and they also developed a series of recovery conditions and error estimates to realize the robust recovery from \eqref{Unconstrained-L1-L2Noise}. However, these results are far from the best. Take for example the mutual coherence condition \footnote{Their original condition in fact is obtained under the cumulative coherence notation. However, whether this condition is sharp was not discussed by the authors.} recently obtained in \cite{li2019signal}, which takes the form of
\begin{align}\label{Li-Chen-MC-condition}
\mu<\frac{2}{\sqrt{3}(5k-2)}.
\end{align}
Obviously \eqref{Li-Chen-MC-condition} is rigorously included in our sharp condition \eqref{Theorem-L2Noise-condition1-Extension-2}. In the aspect of algorithm implementation, since \eqref{unconstrained-L1-DS} is convex, many convex optimization softwares are available to solve it efficiently. Besides, compared to the regularization term (i.e., the second term of the objective function) in \eqref{Unconstrained-L1-L2Noise}, the one in \eqref{unconstrained-L1-DS} is non-smooth and thus non-differentiable. However, if one solves \eqref{unconstrained-L1-DS} using some non-gradient algorithms, such as the alternating
direction method and multipliers \cite{ADMM-2011,wen2017robust}, \eqref{unconstrained-L1-DS} is still comparable to \eqref{Unconstrained-L1-L2Noise} in terms of the algorithmic complexity.

\end{remark}
\begin{IEEEproof}[Proof of Theorem \ref{Theorem-L2Noise-Extension-1}]
The proof is very similar to that of Theorem \ref{Theorem-L2Noise}, and hence we here only present some technique differences. Our proof also relies on Lemma \ref{RIC-like-coherence-Lemma} and the variants of Lemma \ref{NP-Lemma1} and Lemma \ref{NP-Lemma2}. One should keep in mind that the term $\|A\bm{h}\|_{2}$ will be replaced by $\|A^{T}A\bm{h}\|_{\infty}$. Specifically, \eqref{NP-Lemma1-results1}, \eqref{NP-Lemma1-results2} and \eqref{NP-Lemma2-results} are replaced in order by the following inequalities
\begin{align*}
\|A^{T}A\bm{h}\|_{\infty}^{2}-2\epsilon\|A^{T}A\bm{h}\|_{\infty}\leq&2\lambda(\|\bm{h}_{E}\|_{1}-\|\bm{h}_{E^{c}}\|_{1}+2\|\bm{x}_{E^{c}}\|_{1}),\\
\|\bm{h}_{E^{c}}\|_{1}\leq&\|\bm{h}_{E}\|_{1}+2\|\bm{x}_{E^{c}}\|_{1}+\frac{\epsilon}{\lambda}\|A^{T}A\bm{h}\|_{\infty},\\
\|\bm{h}_{E}\|_{2}\leq& \alpha_{1}\|A^{T}A\bm{h}\|_{\infty}+\alpha_{2}\|\bm{h}_{E^{c}}\|_{1}.
\end{align*}
These, as well as the skills in proving Theorem \ref{Theorem-L2Noise}, are sufficient to prove Theorem \ref{Theorem-L2Noise-Extension-1}.
\end{IEEEproof}

\subsection{Structure block-sparse recovery}

Our Theorem \ref{Theorem-L2Noise} can still be extended to guarantee the robust recovery of the structured block-sparse signals. Such a kind of signals (data) arise in many applications \cite{majumdar2010compressed,vidal2006unified,zhang2013compressed}. We assume w.l.o.g. that there are $l$ blocks with block size $d=n/l$ in signal $\bm{x}\in\mathbb{R}^{n}$, and then we can write any signal $\bm{x}\in\mathbb{R}^{n}$ as
\begin{align*}
\bm{x}=[\underbrace{\bm{x}_{1},\cdots,\bm{x}_{d}}\limits_{(\bm{x}[1])^{T}},
\underbrace{\bm{x}_{d+1},\cdots,\bm{x}_{2d}}\limits_{\emph{(\textbf{x}}[2])^{T}},\cdots,
\underbrace{\bm{x}_{N-d+1},\cdots,\bm{x}_{n}}\limits_{\emph{(\textbf{x}}[l])^{T}}]^{T}
\end{align*}
where $\bm{x}[i]\in\mathbb{R}^{d}$ denotes the $i$th block sub-vector of $\bm{x}$. If $\bm{x}$ has at most $k$ non-zero blocks, i.e., $\|\bm{x}\|_{2,0}\leq k$, we refer to such a vector $\bm{x}$ as block $k$-sparse signal. Naturally, a block $k$-sparse signals will reduce the traditional $k$-sparse signal if one takes $d=1$. Accordingly, we can also write any matrix $A\in\mathbb{R}^{m\times n}$ as
\begin{align*}
A=[\underbrace{\bm{a}_{1},\cdots,\bm{a}_{d}}\limits_{A[1]},\underbrace{\bm{a}_{d+1},\cdots,\bm{a}_{2d}}\limits_{A[2]},\cdots,\underbrace{\bm{a}_{N-d+1},\cdots,\bm{a}_{n}}\limits_{A[l]}]
\end{align*}
where $A[i]\in\mathbb{R}^{m\times d}$ denotes the $i$th block sub-matrix of $A$. To recover such a structured block-sparse signal, Eldar amd Mishali \cite{eldar2009robust} proposed solving the following mixed $\ell_{2}/\ell_{1}$-norm minimization problem:
\begin{align}\label{Constrained-L2-L1}
\min_{\bm{x}\in\mathbb{R}^{n}}\|\bm{x}\|_{2,1},~~s.t.~~\bm{b}=A\bm{x},
\end{align}
where $\|\bm{x}\|_{2,1}\triangleq\sum_{i=1}^{l}\|\bm{x}[i]\|_{2}$, and they also derived a block-RIC recovery condition for \eqref{Constrained-L2-L1}. More improved block-RIC conditions can be found in \cite{lin2013block,wang2013recovery,huang2017sharp,li2019high}. As early as 2010, Eldar, Kuppinger and B\"{o}lcskei \cite{eldar2010block} have generalized the traditional mutual coherence to the block setting, and show that any block $k$-sparse signal $\bm{x}$ can be exactly recovered via \eqref{Constrained-L2-L1} if $A$ obeys
\begin{align}\label{BLock-coherence-Condition}
\mu_{B}\leq\frac{1-(d-1)\nu}{(2k-1)d},
\end{align}
where $\mu_{B}$ and $\nu$ are called \emph{block coherence} and \emph{sub-coherence}, respectively, and they are defined as
\begin{align*}
\mu_{B}=\max_{1\leq i<j\leq l}\|(A[i])^{T}A[j]\|_{2} \text{~~and~~} \nu=\max_{1\leq i\leq l}\mu(A[i]).
\end{align*}
Obviously, \eqref{BLock-coherence-Condition} will reduce to \eqref{sharp-mutual-coherence-L1} if one lets the block sub-matrix $A[i]$ be orthonormal \footnote{ This means that $(A[i])^{T}A[i]=I_{d}$, where $I_{d}$ stands for a $d\times d$ identity matrix.} for all $i\in[l]$, namely, $\nu=0$, and also sets the block size $d=1$, see, e.g., \cite{elhamifar2012block,calderbank2015block} for more discussion on block coherence and its related theoretical investigation. Equipped with the block coherence, we here consider extending Theorem \ref{Theorem-L2Noise} to guarantee the robust recovery of such structured block-sparse signals corrupted with the bounded noise by solving the following unconstrained problem
\begin{align}\label{Unconstrained-L1-L2Noise-Block}
\min_{\bm{x}\in\mathbb{R}^{n}}~\|\bm{x}\|_{2,1}+\frac{1}{2\lambda}\|\bm{b}-A\bm{x}\|_{2}^{2}.
\end{align}
This problem sometimes is called the group Lasso \cite{yuan2006model}, and it can also be viewed as the block (group) extension of \eqref{Unconstrained-L1-L2Noise}. One can find our second extension of Theorem \ref{Theorem-L2Noise} as follows.

\begin{theorem}\label{Theorem-L2Noise-Extension-2}
For any $\bm{b}$ observed via \eqref{b=Ax+z} with $\|\bm{z}\|_{2}\leq\epsilon$, if the measurement matrix $A$, whose every block sub-matrix $A[i]$ for $i\in[l]$ is orthonormal, satisfies
\begin{align}\label{Block-coherence-condition-111}
\mu_{\text{B}}<\frac{1}{(2k-1)d}
\end{align}
for certain integer $k\geq2$, then we have
\begin{align*}
\|A(\bm{x}^{\sharp}-\bm{x})\|_{2}&\leq C_{1}(\beta_{1})\|\bm{x}-\bm{x}_{\{k\}}\|_{2,1} + C_{2}(\beta_{1}), \\
\|\bm{x}^{\sharp}-\bm{x}\|_{2}&\leq C_{3}(\beta_{1},\beta_{2})\|\bm{x}-\bm{x}_{\{k\}}\|_{2,1}+ C_{4}(\beta_{1},\beta_{2}),\nonumber
\end{align*}
where $\bm{x}^{\sharp}$ denotes the optimal solution of \eqref{Unconstrained-L1-L2Noise-Block}, $\bm{x}_{\{k\}}$ denotes the best $k$-block approximate of $\bm{x}$, defined as
\begin{align*}
\bm{x}_{\{k\}}=\mathop{\arg\min}_{\|\bm{y}\|_{2,0}\leq k}\|\bm{y}-\bm{x}\|_{2},
\end{align*}
and $\beta_{1}$ and $\beta_{2}$ are defined as
\begin{align*}
\beta_{1}\triangleq\frac{\sqrt{1+(k-1)d\mu_{B}}}{1-(k-1)d\mu_{B}} \text{~~and~~} \beta_{2}\triangleq \frac{\sqrt{k}d\mu_{B}}{1-(k-1)d\mu_{B}}.
\end{align*}
\end{theorem}
\begin{remark}
The idea of using the block coherence and some other tools to deal with the structured block-sparse signals has inspired fruitful results, see, e.g., \cite{lv2011group,wang2014restricted,zeinalkhani2015iterative,hu2017group,zhou2017estimation,elyaderani2017improved,wang2018block}. However, most of these theoretical results focused on the constrained optimization problems rather than their unconstrained counterparts. We note that the authors in \cite{lv2011group} and \cite{elyaderani2017improved} have established some block coherence based theoretical results for an adaptive group Lasso model. Although \eqref{Unconstrained-L1-L2Noise-Block} is included in this adaptive group Lasso model, the results in Theorem \ref{Theorem-L2Noise-Extension-2} are not included in, and are in fact totally different from theirs since Theorem \ref{Theorem-L2Noise-Extension-2} is established in the context of the uniform recovery setting, and the block sparsity requirement of signals is not needed any more when one uses our Theorem \ref{Theorem-L2Noise-Extension-2}, which makes the stable and/or robust recovery of structured block-sparse signals more flexible. Note that one can also extend \eqref{unconstrained-L1-DS} to the block setting and develop a similar theorem with Theorem \ref{Theorem-L2Noise-Extension-2} to deal with the structured block-sparse signals corrupted with the DS type noise.
\end{remark}
\begin{IEEEproof}[Proof of Theorem \ref{Theorem-L2Noise-Extension-2}]
The proof is very similar to that of Theorem \ref{Theorem-L2Noise}, and it relies on the variants of Lemma \ref{RIC-like-coherence-Lemma}, Lemma \ref{NP-Lemma1} and Lemma \ref{NP-Lemma2}. First, \eqref{quasi-RIP-Coherence} will be replaced by
\begin{align*}
(1-(k-1)d\mu_{B})\|\bm{y}\|_{2}^{2}\leq\|A\bm{y}\|_{2}^{2}\leq(1+(k-1)d\mu_{B})\|\bm{y}\|_{2}^{2},
\end{align*}
where $\bm{y}$ represents any block $k$-sparse signal. In fact one can prove it easily using the similar skills in proving Lemma \ref{RIC-like-coherence-Lemma}. Besides, \eqref{NP-Lemma1-results1}, \eqref{NP-Lemma1-results2} and \eqref{NP-Lemma2-results} will also be replaced in order by the following inequalities
\begin{align*}
\|A\bm{h}\|_{2}^{2}-2\epsilon\|A\bm{h}\|_{2}\leq&2\lambda(\|\bm{h}_{E}\|_{2,1}-\|\bm{h}_{E^{c}}\|_{2,1}+2\|\bm{x}_{E^{c}}\|_{2,1}),\\
\|\bm{h}_{E^{c}}\|_{2,1}\leq&\|\bm{h}_{E}\|_{2,1}+2\|\bm{x}_{E^{c}}\|_{2,1}+\frac{\epsilon}{\lambda}\|A\bm{h}\|_{2},\\
\|\bm{h}_{E}\|_{2}\leq& \beta_{1}\|A\bm{h}\|_{2}+\beta_{2}\|\bm{h}_{E^{c}}\|_{2,1},
\end{align*}
where $E$ denotes the block index set over the $k$ blocks with the largest $\ell_{2}$ norm of the original signal $\bm{x}$, and $\beta_{1}$ and $\beta_{2}$ are defined in Theorem \ref{Theorem-L2Noise-Extension-2}. These, as well as the skills in proving Theorem \ref{Theorem-L2Noise}, are sufficient to prove Theorem \ref{Theorem-L2Noise-Extension-2}.
\end{IEEEproof}

\section{Conclusion and Future work}
\label{Conclusions-Futurework}

In this paper, equipped with the powerful mutual coherence tool, we investigated the robust signal recovery using some unconstrained models. We first shown that, if the measurement matrix satisfies $\mu<1/(2k-1)$, one can robustly recover any signal (not necessary to be $k$-sparse) corrupted with the $\ell_{2}$-norm bounded noise using a regularized $\ell_{1}$-norm minimization model \eqref{Unconstrained-L1-L2Noise}. Then we extended this result to guarantee the robust recovery of the signals corrupted with the DS type noise using a DS regularized $\ell_{1}$-norm minimization model \eqref{unconstrained-L1-DS}. To the best of our knowledge, these two kinds of results first extend the sharp uniform recovery condition obtained in \cite{Cai-etal-SharpCoherence-TIT-2010} for \eqref{Noiseless-Constrained-L1} (to guarantee the exact recovery of any $k$-sparse signals) to its two unconstrained variants to guarantee the robust recovery of the signals corrupted with the $\ell_{2}$-norm bounded noise and the DS type noise, respectively. Finally, we considered extending these results to deal with the robust recovery of the structured block-sparse signals corrupted with the bounded noise using some regularized mixed $\ell_{2}/\ell_{1}$-norm minimization models.

There still exists much work to be done in future. Some potential work includes rebuilding the obtained theoretical results using the mutual coherence tool, extending these recovery conditions to guarantee the robust signal recovery in the presence of the random noise, and establishing the coherence-based performance guarantees of some unconstrained convex/nonconvex models for robust vector/matrix/tensor recovery.

%

\ifCLASSOPTIONcaptionsoff
  \newpage
\fi

\end{document}